\documentclass[a4paper,12pt]{amsart}
\usepackage{amsfonts}
\usepackage{mathrsfs}
\usepackage{amsmath,amssymb,latexsym,amsfonts,amscd}
\title{On the geography of threefolds of general type}
\author{Jungkai A. Chen and Christopher D. Hacon}
\address{\rm Taida Institute for Mathematical Sciences,
National Center for Theoretical Sciences, Taipei Office, and
Department of Mathematics, National Taiwan University, Taipei,
106, Taiwan} \email{jkchen@math.ntu.edu.tw}

\address{\rm Department  of Mathematics, University of Utah, 155 South 1400 East, Room 233
Salt Lake City, UT 84112, U.S.A.}
\email{hacon@math.utah.edu}

\thanks{The first author was partially supported by TIMS, NCTS/TPE
and  National Science Council of Taiwan. The second author was
supported by NSF grant
  0456363 and by an AMS Centennial Scholarship.
This work was done during a visit of the first author to the University of Utah.}
\newcommand\Vol{\text{\rm Vol}}

\newcommand\OO{{\mathcal{O}}}

\newtheorem{thm}{Theorem}
\newtheorem{lem}[thm]{Lemma}

\theoremstyle{definition}

\theoremstyle{remark}

\begin{document}
\begin{abstract}
Let $X$ be a complex nonsingular projective 3-fold of general type.
We show that there are positive constants $c$, $c'$ and $m_1$ such that $\chi (\omega _X)\geq -c\Vol (X)$ and $P_m(X)\geq c'm^3\Vol (X)$ for all $m\geq m_1$.
\end{abstract}
%%%%%%%%%%%%%%%%%%%%
\maketitle
%%%%%%%%%
\pagestyle{myheadings} \markboth{\hfill J. A. Chen and C. D. Hacon \hfill}{\hfill On the geography of threefolds of general
type\hfill}
%%%%%%%%%%%%%%%%%%%%%%%%%%%%%%%%%%

\section{Introduction and known results}
The birational classification of surfaces of general type is well understood.
For example, it is known that if $X$ is a surface of general type then $|mK_X|$ induces a birational map for all $m\geq 5$.
As a general rule, it is not possible to classify surfaces of general type with given invariants.
%There are of course notable exception such as Castelnuovo's Criterion according to which if $X$ is a surface with $P_2(X)=q(X)=0$, then $X$ is birational to$\mathbb P ^2$.
In general, the best that one can do is to show that
the invariants of a surface $X$ satisfy certain
inequalities.
A fundamental inequality for the invariants of a minimal surface of general type is the Bogomolov-Miyaoka-Yau inequality $K_X^2\leq 9\chi (\OO _X)$.
%Two fundamental inequalities for the invariants of a minimal surface of general type are: the Bogomolov-Miyaoka-Yau inequality $K_X^2\leq 9\chi (\OO _X)$ and the  Noether inequality $K_X^2\geq 2p_g-4$.

It is a natural problem to try and extend the results for surfaces to higher dimensions. There have been many partial results for $3$-folds.
For example, it is shown in \cite{CCM} that if $X$ is a Gorenstein minimal $3$-fold of general type, then $|mK_X|$ induces a birational map for all $m\geq 5$.
In fact the proof is based upon the fact that for such $3$-folds, we have
the Miyaoka-Yau inequality $K_X^3\leq 72\chi (\omega _X)$.% or equivalently $K_X^3\leq  3K_X\cdot c_2(X)$. %In \cite{CCZ}, it is also shown that the Noether type inequality $K_X^3\geq \frac 23(2p_g-5)$ also holds (for smooth minimal $3$-folds of general type). These inequalities are very useful analogoues to the Bogomolov-Miyaoka-Yau and Noether inequalities for surfaces. The main drawback is that they do not hold for $3$-folds of index $>1$.

Despite many partial results, the geometry of non-Gorenstein 3-folds of general type has proven to be a very challenging topic.
In a recent paper however, the first author and M. Chen \cite{CC} show the remarkable result that if $X$ is a smooth complex projective $3$-fold of general type,
then $P_{12}\geq 1$, $P_{24}\geq 2$ and $|mK_X|$ induces a birational map for all $m\geq 77$. It is then natural to hope that further precise results
on the geography of $3$-folds of general type may be within reach.
The purpose of this paper is to show that using the methods of \cite{CC}
one can in fact prove an inequality similar to the Miyaoka-Yau inequality
which holds for non-Gorenstein $3$-folds of general type. Namely we show that
\begin{thm}\label{T-1} There exists a constant $c>0$ such that for any minimal $3$-fold of general type with terminal singularities $X$, we have
$$\chi(\omega _X)\geq -c K_X^3.$$\end{thm}
Recall that for any minimal $3$-fold of general type with terminal singularities, we have ${\rm Vol}(X)=K_X^3$.
It should be noted that $\chi(\omega _X)$ may be negative
for $3$-folds of general type.
In fact consider curves $C_1$, $C_2$ and $C_3$ with genus $g_i$ and
involutions $\sigma _i$ such that $C_i/<\sigma _i >\cong \mathbb P ^1$.
Then the $3$-fold $X$ given by a desingularization of the quotient of
$C_1\times C_2\times C_3$ by the ``diagonal'' involution, has $\chi (\omega _X)<0$.
In fact if we let $g_1=g_2=g$, then for fixed $g_3$ and for $g\gg 0$ one has that
$-\chi (\omega _X)=O(g^2)$ and $K_X^3=O(g^2)$. So the inequality of Theorem
\ref{T-1} has the right shape. The constant $c$
that may be computed with the methods of this paper and the results of
\cite{CC} is $c=32\cdot 120^3$. We expect that this is far from optimal
and so we make no effort to determine it explicitly.
We remark that if $\Vol(X)\gg 0$, then using the results of \cite{T}, one
can recover $c=2502$.

We also prove the following result concerning the plurigenera of $X$.
\begin{thm}\label{T-2} There exist constants $c'>0$ and $m_1>0$ such that for any minimal $3$-fold of general type with terminal singularities $X$, we have $$P_m(X)\geq c'm^3K_X^3\qquad {\rm for\  all}\ m\geq m_1.$$
\end{thm}
Once again the values of $c'$ and $m_1$ that may be computed with the methods of this paper are far from optimal so we make no effort to determine their values (note that using the results of \cite{CC}, it follows form the arguments below that $c'=\frac 5{89168}$ and $m_1=112$ suffice).

It would be interesting to: \begin{enumerate}
\item Determine the optimal value of $m_1$ in Theorem \ref{T-2};
\item See if Theorem \ref{T-2} can be recovered by using
the methods of \cite{T};
\item See if Theorems \ref{T-1} and \ref{T-2} hold in higher dimensions.\end{enumerate}

We remark that the proof of the results of this paper is based on the methods of \cite{CC}.
We have chosen to keep the exposition of this paper as self contained and simple as possible. Therefore, we include a proof of all the results
of \cite{CC} that we will use (namely inequalities (1) and (2)).
\section{Some inequalities}
In this section, we will prove the following inequalities: \begin{thm}Let $X$ be a minimal $3$-fold of general type with terminal singularities, then $$P_4+P_5+P_6-3P_2-P_3-P_7 \ge 0.
\eqno(1)$$ and
$$ 2P_5+3P_6+P_8+P_{10}+P_{12} \ge \chi(\OO _X) +10P_2+4 P_3
+P_7+P_{11}+P_{13}+14\sigma _{12} \eqno(2)$$
where $\sigma _{12}$ is a positive integer that will be defined below.\end{thm}
Note that a stronger version of the above inequalities is proved in \cite{CC}.
Here we include a simpler and self contained proof of this (weaker) version of the inequalities of \cite{CC}. The stronger version also follows from
the methods of this paper, however it is not necessary for our purposes so we
have chosen not to include it here.

We consider now $X$ a minimal $3$-fold of general type with
terminal singularities. According to Reid (see last section of
\cite{YPG}), there is a  ``basket'' of pairs of integers
$\mathscr{B}(X):=\{(b_i, r_i)\}$ such that the Riemann-Roch formula may be written as
$$\chi(\mathcal{O}_X(mK_X))=\frac{1}{12}m(m-1)(2m-1)K_X^3-(2m-1)\chi(\mathcal
{O}_X)+l(m),  $$ where the correction term $l(m)$ is computed
by:
$$l(m):=\sum_{Q_i\in \mathscr{B}(X)}l_{Q_i}(m):=\sum_{Q_i\in \mathscr{B}(X)}\sum_{j=1}^{m-1}
\frac{\overline{jb_i}(r_i-\overline{jb_i})}{2r_i}.$$ Here, we assume that
$b_i$ is
co-prime to $r_i$ and $0<b_i\leq \frac{r_i}{2}$.  The ratio $\frac{b_i}{r_i}$ is
called the {\it slope of $(b_i, r_i)$}.
For a basket $B$, we let $\sigma (B):=\sum b_i$ and $\sigma_{12} (B):=\sum_{\frac{b_i}{r_i}\leq \frac 1 {12}} b_i$.

Let
$$\overline{M}^j(b,r):=\frac{\overline{jb}(r-\overline{jb})}{2r}, \quad
{M^j}(b,r):=\frac{{jb}(r-{jb})}{2r},$$
 $$\Delta^j({b,r}):=\overline{M}^j(b,r)-M^j(b,r).$$ An easy computation
shows that $\Delta^n({b,r})=ibn-\frac{i^2+i}{2}r$, where $i=\lfloor
\frac{bn}{r} \rfloor$.

We will need the following easy computational lemmas.

\begin{lem} \label{nodiff} Let $b_1r_2-b_2r_1=1$. If $0<n\ne xr_1+yr_2$
for any integers $x,y >0$, then there is no rational number $\frac{b}{n} \in (\frac{b_2}{r_2},\frac{b_1}{r_1})$ and we have
$$\Delta^n({b_1+b_2,r_1+r_2})=\Delta^n({b_1,r_1})+\Delta^n({b_2,r_2}).$$
\end{lem}
\begin{proof} We may assume that $\frac {b_2}{r_2}<\frac{b_1}{r_1}$.
Note also that by our assumptions, we have $n\leq r_1r_2$.
%Note that since $r>1$, we have $\frac{b_1n}{r_1}-\frac {b_2n}{r_2}=\frac n{r_1r_2}$.
If $\frac{b}{n} \in (\frac{b_2}{r_2},\frac{b_1}{r_1})$, then $n= (br_2-b_2n)r_1+(b_1n-br_1)r_2$ with  $br_2-b_2n>0$ and $b_1n-br_1 >0$.
Hence, as $n \ne xr_1+yr_2$ for all $x,y >0$, then  there is no rational number $\frac{b}{n} \in (\frac{b_2}{r_2},\frac{b_1}{r_1})$.

Let $i_1:=\lfloor
\frac{b_1n}{r_1} \rfloor$, $i_2:=\lfloor
\frac{b_2n}{r_2} \rfloor$ and $i:=\lfloor
\frac{(b_1+b_2)n}{r_1+r_2} \rfloor$. If $\frac{b_1n}{r_1}$ is not an integer, then $$i_2=\lfloor
\frac{b_2n}{r_2} \rfloor \le i_1=\lfloor
\frac{b_1n}{r_1} \rfloor<\frac{b_1n}{r_1}.$$ If $i_1\ne i_2$
then $i_1>\frac {b_2n}{r_2}$ so that $\frac{i_1}n\in (\frac{b_2}{r_2},\frac{b_1}{r_1})$ which is impossible. Therefore $i_1=i_2=i$ and the statement follows from the equation $\Delta=ibn-\frac{i^2+i}2r$.

If $\frac{b_1n}{r_1}$ is an integer, then one sees that $\Delta ^n(b_1,r_1)
=(i_1-1)b_1n-\frac{(i_1-1)^2+(i_1-1)}2r_1$ so that as $i_2=i_1-1$, the statement follows from the definition of $\Delta$.
\end{proof}

\begin{lem} \label{diff} Suppose that $b_1r_2-b_2r_1=1$ and
suppose that $n=xr_1+yr_2$ for some integers $r_2 \ge x >0$, $r_1 \ge
y >0$, then
$$\Delta^n({b_1+b_2,r_1+r_2})=\Delta^n({b_1,r_1})+\Delta^n({b_2,r_2})-\min\{x,y\},$$
\end{lem}

\begin{proof}
We first remark that the expression $n=xr_1+yr_2$ for some integers $r_2 \ge x >0$, $r_1 \ge
y >0$ is unique.

Let $i=xb_1+yb_2$.\\
An easy computation shows that if $r_1\ne y$, then
$$\lfloor \frac{b_1n}{r_1} \rfloor = \lfloor \frac{xb_1r_1+yb_1r_2 }{r_1} \rfloor =\lfloor \frac{xb_1r_1+yb_2r_1+y }{r_1} \rfloor =i+\lfloor \frac{y }{r_1} \rfloor=i,$$

$$\lfloor \frac{b_2n}{r_2} \rfloor = \lfloor \frac{xb_2r_1+yb_2r_2 }{r_2} \rfloor =\lfloor \frac{xb_1r_2-x+yb_2r_2 }{r_2} \rfloor =i+\lfloor \frac{-x }{r_2} \rfloor=i-1,$$

$$\lfloor \frac{(b_1+b_2)n}{r_1+r_2} \rfloor = \lfloor \frac{xb_1r_1+yb_1r_2+xb_2r_1+yb_2r_2 }{r_1+r_2} \rfloor $$
$$=\lfloor \frac{xb_1r_1+yb_2r_1+y+xb_1r_2-x+yb_2r_2 }{r_1+r_2} \rfloor =i+\lfloor \frac{y-x }{r_1+r_2} \rfloor.$$

If $y \ge x$, then $\lfloor \frac{(b_1+b_2)n}{r_1+r_2} \rfloor=i$. Direct computation gives
$$ \Delta^n({b_1+b_2,r_1+r_2})-\Delta^n({b_1,r_1})-\Delta^n({b_2,r_2}) $$
$$=i(b_1+b_2)n-\frac{i^2+i}{2}(r_1+r_2)-ib_1n + \frac{i^2+i}{2}r_1 -(i-1)b_2n+\frac{i^2-i}{2}r_2$$
$$ = b_2n-ir_2 = b_2(xr_1+yr_2)-(xb_1+yb_2)r_2=x(b_2r_1-b_1r_2)=-x.$$
Note that if $r_1=y$, then $\lfloor \frac{b_1n}{r_1} \rfloor =i+1$.
However, one easily sees that
the above formula is unchanged.

If $y \le x$, the computation is similar.
\end{proof}

\begin{proof}[Proof of inequality (1)]
By direct computation, one finds that the $K_X^3$ and $\chi(\OO _X)$ terms coming from the Riemann-Roch formula cancel.
Inequality (1) is then equivalent to
$$ -3 l(2)-l(3)+l(4)+l(5)+l(6)-l(7) \ge 0.$$
Since $l(m)=\sum_{j=1}^{m-1} \overline{M}^j(B)=\sum_{j=1}^{m} \sum
\overline{M}^j(b_i,r_i)$, we must show the inequality
$$
\overline{\Xi}(B):=-2\overline{M}^1(B)+\overline{M}^2(B)+2\overline{M}^3(B)+\overline{M}^4(B)-\overline{M}^6(B)
\ge 0. \eqno(3)$$

We will show that this holds for any single basket $(b,r)$
and hence for any basket
$B$.

We define $$\Xi(B):=-2{M^1}(B)+{M^2}(B)+2{M^3}(B)+{M^4}(B)-{M^6}(B),$$
$$\Xi\Delta(B):=\overline{\Xi}(B)-\Xi(B)=2{\Delta^3}(B)+{\Delta^4}(B)-{\Delta^6}(B),$$
where we have used the fact that as we assumed that $b/r\leq 1/2$, then
${\Delta^1}(B)={\Delta^2}(B)=0$.
%Clearly, $\Xi\Delta(B)+\Xi(B)=\overline{\Xi}(B)$.

 {\bf Step 1.} For any single basket $B=\{(b,r)\}$, we have
$\Xi(B)=2b$.

{\bf Step 2.} For the single basket $B=\{(1,2)\}$, we have $\Delta^3({1,2})=1,\Delta^4({1,2})=2$ and $\Delta^6({1,2})=6$.
Hence
$\Xi\Delta(B)=-2$, and $\overline{\Xi}(B)=0$. A similar computation for
$B=\{(1,3)\}$ or $B=\{(1,4)\}$, yields $\Xi\Delta(B)=-2$ and $\overline{\Xi}(B)=0$. When $B=\{(1,5) \}$, then $\Xi\Delta(B)=-1$
and $\overline{\Xi}(B)=1$.

{\bf Step 3.} When $B=\{(1,r)\}$ with $r \ge 6$, we have $\overline{M}^m(1,r)=M^m(1,r)$ for all $m \le 6$. Hence $\overline \Xi(B)=\Xi(B)=2$.

{\bf Step 4.}
Recall that we are assuming $\frac{b}{r} \le \frac{1}{2}$.
Let $S=\{\frac{1}{r}\}_{r \ge 2}$, $S^{(5)}:=S \cup \{\frac{2}{5}\}$ and
for $n\geq 6$ set
$$S^{(n)}=S^{(n-1)} \cup \{\frac{b}{n}| (b,n)=1,0< \frac{b}{n} \le \frac{1}{2} \}.$$ For any $\frac{b}{n} \in S^{(n)}$, let $[0;a_1,...,a_t]$ be its continued fraction expression. Note that as $\frac bn\leq \frac 12$, then $a_1\geq 2$.
If $t>1$, we may consider the rational number
$\frac{b_1}{r_1}$ with continued fraction expression $[0;a_1,...,a_{t-1}]$.
We have that $nb_1-r_1b = \pm 1$ and $\frac{b_1}{r_1}\leq \frac 12$.
Let $b_2=b-b_1, r_2=n-r_1$. Notice that we also
have $b_1r_2-b_2r_1=\pm 1$ and $\frac{b_2}{r_2}\leq \frac 12$.
Then we have $\frac{b_1}{r_1}, \frac{b_2}{r_2} \in S^{(n-1)}$.

{\bf Step 5.} We proceed by showing by induction on $r$ that inequality (3) holds. By Step 1, this is equivalent to showing that
$\Xi \Delta (b,r)\geq -2b$.

We have seen that the inequality (3) holds for $r \le 4$. For $r=5$, we must
consider the single basket $B=\{(2,5)\}$. Notice that
$\Delta^n({2,5})=\Delta^n(1,2)+\Delta^n({1,3})$, for $n=3,4,6$ by Lemma \ref{nodiff}. We see that
$$\Xi\Delta(2,5)=\Xi\Delta(1,2)+\Xi\Delta(1,3)=-4.$$
By Step 1, we have $\overline{\Xi}(2,5)=0$.

For $r=6$, there are no new baskets to consider.

{\bf Step 6.}
For $r \ge 7$, notice that by Step 4, we may assume that $(b,r)=(b_1,r_1)+(b_2,r_2)$ for some
$r_1,r_2 < r$ and (after possibly switching indices) that $b_1r_2-b_2r_1=1$. By induction hypothesis, we have
$\Xi\Delta(b_i,r_i) \ge -2 b_i$. Using Lemma \ref{nodiff}, it is easy to see that
$$\Delta^m(b,r)=\Delta^m(b_1,r_1)+\Delta^m(b_2,r_2)$$for $m \in \{ 3,4,6\}$.
Hence $$\Xi\Delta(b,r)=\Xi\Delta(b_1,r_1)+\Xi\Delta(b_2,r_2) \ge -2b.$$

This completes the proof.
\end{proof}

\begin{proof}[Proof of inequality (2)]
The proof is similar but the computations are a little bit more involved.

Inequality (2) is  equivalent to
$$ -10 l(2)-4 l(3)+2 l(5)+3 l(6)-l(7)+l(8)+l(10)-l(11)+l(12)-l(13) \ge 14\sigma _{12},$$
which in turn is equivalent to
$$
\overline{\Xi}(B) := -9\overline{M}^1(B)+\overline{M}^2(B)+5\overline{M}^3(B)+5\overline{M}^4(B)$$
$$\qquad +3\overline{M}^5(B)
+\overline{M}^7(B)-\overline{M}^{10}(B)-\overline{M}^{12}(B) \ge 14\sigma _{12}(B). \eqno(4)$$

We will show that this holds for any single basket and hence for any
basket $B$.

We define $\Xi(B)$ and $\Xi\Delta(B)$ as in the proof of inequality (1).

 {\bf Step 1.} For any single basket $B=\{(b,r)\}$, we have
$\Xi(B)=14b$.

{\bf Step 2.} For a single basket $B=\{(1,r)\}$ with $2\leq r \le 11$,
direct computation gives $\overline{\Xi}(1,r)=0,0,0,2,5,6,8,10,12,13$.

{\bf Step 3.} We claim that if $B=\{(b,r)\}$ with $\frac b r\leq \frac 1{12}$, then
$\overline \Xi(B)=14b$.

When $B=\{(1,r)\}$ with $r \ge 12$, we have $\overline{M}^m(1,r)=M^m(1,r)$ for all $m \le 12$, therefore $\overline \Xi(B)=\Xi(B)=14$.
When $B=\{(b,r)\}$ with $\frac b r<\frac 1{12}$ and $b>1$, as in the proof of
inequality (1) Step 4, we
may write $b=b_1+b_2$ and $r=r_1+r_2$ where $b_i$ and $r_i$ are co-prime,
$\frac{b_1}{r_1}>\frac b r > \frac{b_2}{r_2}$ and $b_1r_2-b_2r_1=1$.
By Lemma \ref{nodiff}, we see that as $r=r_1+r_2>12$, then
$\frac 1 {12}\not\in(\frac{b_2}{r_2},\frac{b_1}{r_1})$.
It follows that $\frac{b_1}{r_1}\leq \frac 1{12}$. The claim now follows by induction. In fact, since $r>12$, by Lemma \ref{nodiff}, we have
$$\Delta ^m(b,r)=\Delta ^m(b_1,r_1)+\Delta ^m(b_2,r_2)$$ for all $1\leq m\leq 12$.

We proceed by showing by induction on $r$ that $\Xi \Delta (b,r)\geq -14b$. By Step 1, this is equivalent to $\overline{\Xi} (b,r)\geq 0$ and hence implies that inequality (4) holds.

{\bf Step 4.} $\overline{\Xi} (b,r)\geq 0$ for all baskets $(b,r)$ with $r\leq 12$.

By Step 1, $\overline{\Xi} (b,r)\geq 0$ for all baskets $(b,r)$ with $r \le 4$. For $r=5$, we must
consider the single basket $B=\{(2,5)\}$. By Lemmas \ref{nodiff} and \ref{diff},
one sees that
$$\Delta^n({2,5})-\Delta^n(1,2)-\Delta^n({1,3})=-1,-1,-2,-2\quad {\rm for}\ n=5,7,10,12$$ respectively and $\Delta^n({2,5})-\Delta^n(1,2)-\Delta^n({1,3})=0$
for $n=1,2,3,4$. It follows that
$$\Xi\Delta(2,5)=\Xi\Delta(1,2)+\Xi\Delta(1,3).$$
By Steps 1 and 2, we have $\overline{\Xi}(2,5)=0$.

For $r=6$, there are no new baskets to consider.

We can similarly compute all single baskets $B\in S^{(12)}-S^{(6)}$. Recall that each single basket $(b,r)$,
can be compared with pairs $(b_1,r_1)$ and $(b_2,r_2)$ as described in Step 4 of the proof of inequality (1).

We have that $\Xi\Delta(b,r)\geq \Xi\Delta(b_1,r_1)+\Xi\Delta(b_2,r_2)$ for all
$B\in S^{(12)}-S^{(6)}$ or more precisely that
$\Xi\Delta(b,r)=\Xi\Delta(b_1,r_1)+\Xi\Delta(b_2,r_2) +1$ if $(b,r)\in \{(3,10),(5,12)\}$ and $\Xi\Delta(b,r)=\Xi\Delta(b_1,r_1)+\Xi\Delta(b_2,r_2)$
otherwise. Therefore $\Xi\Delta(b,r)\geq -14b$ for all baskets $(b,r)$ with $r\leq 12$.

{\bf Step 5.} $\overline{\Xi} (b,r)\geq 0$ for all baskets $(b,r)$ with
$r \ge 13$ and $\frac b r >\frac 1{12}$.

We may assume that $(b,r)=(b_1,r_1)+(b_2,r_2)$ for some
$r_1,r_2 < r$. By induction hypothesis, we have
$\Xi\Delta(b_i,r_i) \ge -14 b_i$. By Lemma \ref{nodiff}, we have
$$\Delta^m(b,r)=\Delta^m(b_1,r_1)+\Delta^m(b_2,r_2),$$
for $m \le 12$. Hence $$\Xi\Delta(b,r)=\Xi\Delta(b_1,r_1)+\Xi\Delta(b_2,r_2) \ge -14b.$$

This completes the proof.
\end{proof}
We will also need the following equality
\begin{lem}\label{l-6} For any minimal $3$-fold of general type with terminal
singularities and basket $B$ we have $\sigma (B)=10\chi (\OO _X)+5P_2(X)-P_3(X)$.
\end{lem}
\begin{proof} The equality follows immediately from the Riemann-Roch formula.
\end{proof}

\section{Main result}
In this section we prove the main result of this paper.
\begin{thm}
Let $X$ be a smooth $3$-fold of general type. Then
\begin{enumerate}
\item There are constants $c'>0$ and $m_1>0$ such that $P_m(X) \ge c' m^3\Vol(X)$ for all $m \ge m_1$.

\item There is a constant $c>0$ such that $\Vol(X) \ge c \chi(\OO_X)$.
\end{enumerate}
\end{thm}
\begin{proof}
We will first prove (1). Consider the Riemann-Roch formula.

If $\chi(\OO_X) \le 0$. Then we get
$$P_m \ge \frac{m(m-1)(2m-1)}{12} \Vol(X) \ge \frac{m^3}{16}\Vol(X)$$ for $m \ge 2$ already.

It remains to consider the case when $\chi(\OO_X) >0$. We will need the following.

\begin{lem} There exist constants $m_0,c_1,c_2 >0$ such that
\begin{enumerate} \item $P_{m_0}\geq 2$,
\item  $P_m\geq 2$ for all $m\geq 5m_0+6$,
\item $P_m \ge c_1{m}$ for all $m\geq 12m_0+10$
and \item $P_m \ge \frac{c_2 m}{t} P_t$ for any  $m \ge 10m_0+2t+10$.\end{enumerate}\end{lem}
\begin{proof} We will repeatedly use the fact that if $P_s>0$ and $P_t>0$, then $P_{s+t}\geq P_s+P_t-1$ and so for all $s\geq t_0=5m_0+6$ and any $t'>0$ such that $P_{t'}\geq 2$, we have
$$P_s> \lfloor \frac{s-t_0}{t'}\rfloor (P_{t'}-1)\geq \frac{s-t_0-t'+1}
{t'}(P_{t'}-1).$$

(1) If $P_i\leq 1$ for $i\in \{ 5,6,8,10,12 \}$, then by inequality (2), we have
$0< \chi(\OO _X)\leq 8$ and $\sigma _{12}=0$.
Since $\sigma =10\chi(\OO _X)+5P_2-P_3$ (cf. Lemma \ref{l-6}), we have
$\sigma=\sum b_i \leq 85$.
Therefore as $\sigma _{12}=\sum _{\frac {b_i}{r_i}\leq \frac 1 {12}}b_i=0$,
there are only finitely many possible
such baskets of singularities and hence there is an integer $m_0$ such that $P_{m_0}(X)\geq 2$.
We may assume that $120$ divides $m_0$.

(2) If $P_i\ge 2$ for some $i\in \{ 5,6,8,10,12 \}$, then we have $P_{120}(X)\geq 2$.
Therefore, if $\chi (\OO _X)>0$, then
$P_{m_0}(X)\geq 2$. By \cite{MC} we have that $|mK_X|$ is birational for all $m\geq 5m_0+6$.

(3) It follows that for all $m\geq 12m_0+10$ we have
$$P_m> \frac{m-6m_0-5}{m_0}(P_{m_0} -1)\geq \frac {m}{2m_0}.$$

(4) If $P_t=0$, the proposed inequality is trivial. If $P_t=1$, the proposed inequality follows from (3) assuming that $c_2\leq c_1$. We now assume that
$P_t\geq 2$ and hence $P_t-1\geq P_t/2$. We have that for $m\geq 10m_0+2t+10$,
$$P_m> \frac{m-5m_0-t-5}{t}(P_t -1)\geq \frac {m} {4t} P_t.$$
\end{proof}

We have that $2P_5+3P_6+P_8+P_{10}+P_{12} \ge \chi(\OO_X)$.
Hence $$(2\frac{5}{c_2 m} + 3 \frac{6}{c_2 m} +\frac{8}{c_2 m}+\frac{10}{c_2 m} + \frac{12}{c_2 m})P_m=\frac{58}{c_2m} P_m\ge \chi(\OO_X)$$
for any $m\geq 10m_0+34$.

Thus, by the Riemann Roch formula $$ (1+\frac{116}{c_2})P_m \ge   P_m+2m \chi(\OO_X) \ge \frac{m^3}{16}\Vol(X).$$
This proves the first inequality.

The second inequality holds trivially if $\chi(\OO_X) \le 0$. Hence we assume that $\chi(\OO_X) >0$.
If $P_5,P_6,P_8,P_{10},P_{12} \le 1$, then $\chi(\OO_X) \le 8$.
Since $|(5m_0+6)K_X|$ is birational, then $\Vol (X)\geq \frac 1 {(5m_0+6)^3}$.
Therefore,
$$ \Vol(X) \ge \frac{1}{ (5m_0+6)^3    } \ge \frac{1}{(5m_0+6)^3      \cdot 8} \chi(\OO_X).$$
%Therefore, by \cite[Theorem 1.2]{CC},$$ \Vol(X) \ge \frac{1}{2660} \ge \frac{1}{2660 \cdot 8} \chi(\OO_X).$$ Note that instead of the inequality  $\Vol(X) \ge \frac{1}{2660}$, one could use the fact that as $|(5m_0+6)K_X|$ is birational, then $Vol(X)\geq \frac 1 {5m_0+6)^3}$.

In general, we have $P_{120} \ge P_t$ for $t\in \{5,6,8,10,12\}$.
Hence $8P_{120} \ge \chi(\OO_X)$. We may assume that $P_{t}\geq 2$ for some $t\in \{5,6,8,10,12\}$ so that
$|120K_X|$ is birational. Therefore
$120^3\Vol(X) \ge P_{120}-3\geq 1$, hence $$4\cdot 120^3\Vol(X) \ge 120^3\Vol(X)+3 \ge  P_{120} \ge \frac{1}{8} \chi(\OO_X).$$
\end{proof}
%%%%%%%%%%%%%%%%%%%%%%%%%%%%%%%%%%%%%%%%%%%%%%%%%%%%%%%

\end{document}